Research Paper

# Inter-laboratory replicability and sensitivity study of a finite element model to quantify human femur failure load: case of metastases

Marc Gardegaront [1, 2], Amelie Sas [3], Denis Brizard [2], Aurélie Levillain [2], François Bermond [2], Cyrille B. Confavreux [1,4], Jean-Baptiste Pialat [4, 5], G. Harry van Lenthe [3], Hélène Follet [1] and David Mitton* [2]

[1] Univ Lyon, Univ Claude Bernard Lyon 1, INSERM, LYOS UMR 1033, 69008 Lyon, France; m.gardegaront@gmail.com; helene.follet@inserm.fr

[2] Univ Lyon, Univ Eiffel, Univ Claude Bernard Lyon 1, LBMC UMR_T9406, 69622 Lyon, France; denis.brizard@univ-eiffel.fr; aurelie.levillain@univ-lyon1.fr; francois.bermond@univ-eiffel.fr; david.mitton@univ-eiffel.fr

[3] Biomechanics Section, Dept. Mechanical Engineering, KU Leuven, Leuven, Belgium; amelie.sas@kuleuven.be; harry.vanlenthe@kuleuven.be

[4] Centre Expert des Métastases Osseuses (CEMOS), Hôpital Lyon Sud, Hospices Civils de Lyon ; cyrille.confavreux@chu-lyon.fr; jean-baptiste.pialat@chu-lyon.fr

[5] Creatis CNRS UMR 5220, INSERM U1294, Université Lyon 1, Villeurbanne, France

* Correspondence: david.mitton@univ-eiffel.fr

**Highlights**

Reproducibility of the femur failure load model is validated

Replicability of the femur model shows different performances on two datasets

Sensitivity of the model with respect to density-based parameters is high

Sensitivity with respect to segmentation, orientation and femur length is not negligible

**Abstract**

Introduction: Metastases increase the risk of fracture when affecting the femur. Consequently, clinicians need to know if the patient's femur can withstand the stress of daily activities. The current tools used in clinics are not sufficiently precise. A new method, the CT-scan-based finite element analysis, gives good predictive results. However, none of the existing models were tested for reproducibility. This is a critical issue to address in order to apply the technique on a large cohort around the world to help evaluate bone metastatic fracture risk in patients. The aim of this study is then to evaluate 1) the reproducibility 2) the replicability and 3) the global sensitivity of one of the most promising models of the literature (original model). Methods: The reproducibility was evaluated by comparing the results given in the original model by the original first team (Leuven, Belgium) and the reproduced model made by another team (Lyon, France) on the same dataset of CT-scans of *ex vivo* femurs. The replicability was evaluated by comparing the results of the reproduced model on two different datasets. The global sensitivity analysis was done by using the Morris method and evaluates the influence of the density calibration coefficient, the segmentation, the orientations and the length of the femur. Results: The original and reproduced models are highly correlated ($r^2 = 0.95$), even though the reproduced model gives systematically higher failure loads. When using the reproduced model on another dataset, predictions are less accurate ($r^2$ with the experimental failure load decreases, errors increase). The global sensitivity analysis showed high influence of the density calibration coefficient (mean variation of failure load of 84 %) and non-negligible influence of the segmentation, orientation and length of the femur (mean variation of failure load between 7 and 10 %). Conclusion: This study showed that, although being validated, the reproduced model underperformed when using another dataset. The difference in performance depending on the dataset is commonly the cause of overfitting when creating the model. However, the dataset used in the original paper [18] and the Leuven's dataset gave similar performance, which indicates a lesser probability for the overfitting cause. Also, the model is highly sensitive to density parameters and automation of measurement may minimize the uncertainty on failure load. An uncertainty propagation analysis would give the actual precision of such model and improve our understanding of its behavior and is part of future work.

**Graphical abstract**

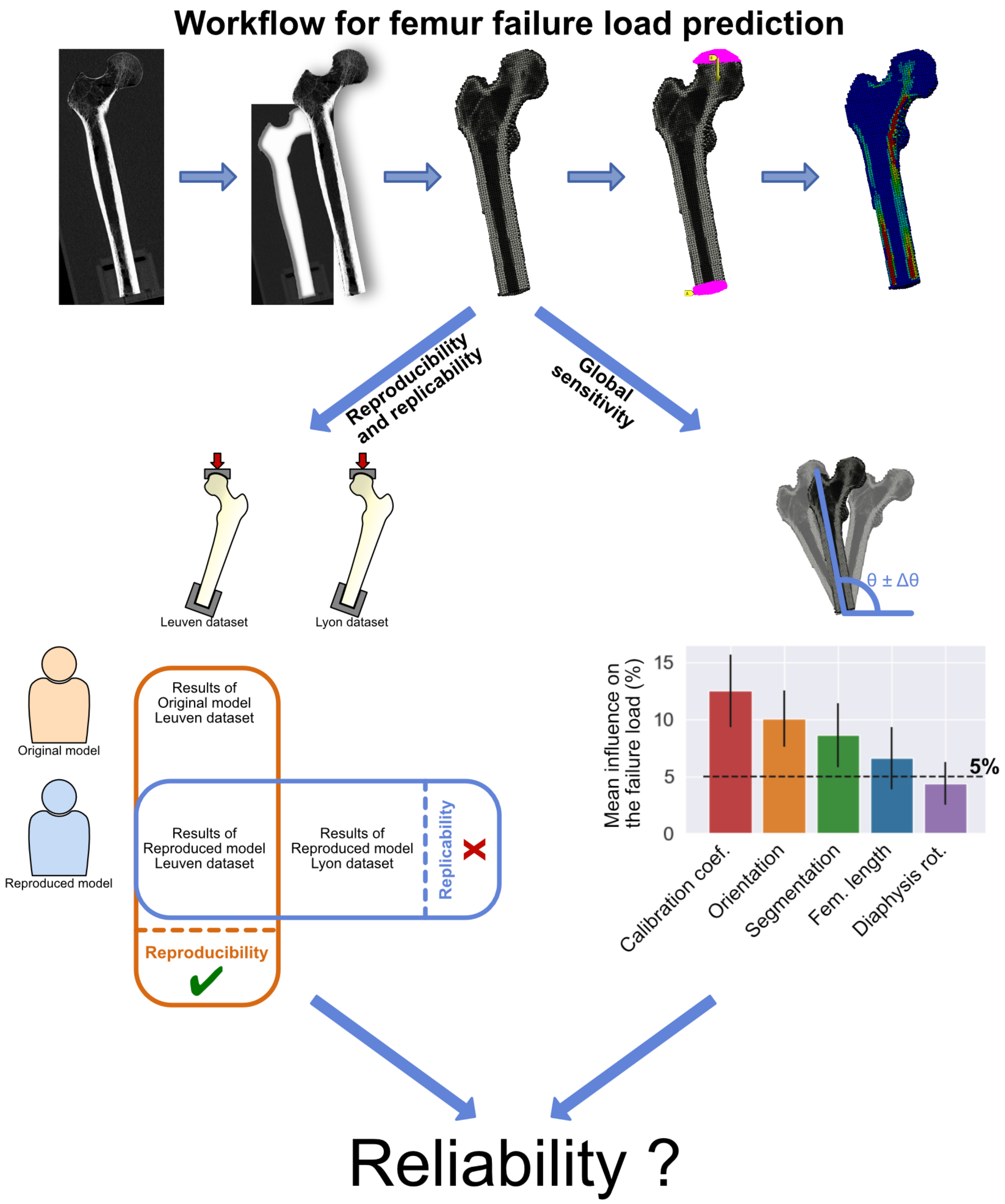

## Introduction

Bone metastases occur when cancer cells from a primary tumor spread to the bone through the bloodstream or lymphatic system. Bone is an important site for metastases [1]. When affecting weight-bearing bones such as femurs, the metastasis increases the risk of fracture occurring during normal activities. A fracture occurring on already enfeebled cancer patients may lead to complications that will impair their oncologic treatments, hinder their chances of survival, and decrease their quality of life [2]. Assessing the risk of fracture is critical to decide whether the femur of a specific patient can handle daily activities or if it has to undergo prophylactic surgery to avoid fracture complications. This surgery, besides being less harmful than an actual fracture, impacts nonetheless the quality of life of the patient for several weeks/months and delays oncological treatment. Therefore, the objective for clinicians is to deliver prophylactic surgery at the appropriate time.

Clinicians have developed tools to help them predict impending fractures [1]. The current gold standard for metastasized femur fracture prediction is the Mirels' score [3]. However, several researchers have shown that this score is not sufficiently precise [4 - 7], and its strict application could lead to overestimation of the risk of fracture. Therefore, other models have been created to better predict impending fracture. Among them, computed tomographic (CT) scan-based finite element (FE) models have given good predictive results from various teams since 1998 and have improved further since then [8 - 18]. These FE models were all validated by their respective authors, either by comparing the predictions to *ex vivo* experiments or, less commonly, by confronting the model predictions to clinicians' predictions [16]. However, to our knowledge, none of these models were evaluated for replicability.

Many years of research have established that reproducibility (obtaining the same results based on the same data and methods by two independent investigators) and replicability (obtaining the same performance based on different data) [19] are challenging topics in science [20 - 22], and too few teams publish analyses of global sensitivity of their FE model and disclose the values of their FE model results and parameter uncertainties. To our knowledge, only two papers explicitly discussed partial model sensitivity to some parameters [23] [24] and none examined their replicability. It has already been established that results may depend on the operator [25] without even replicating it.

Therefore, as a first step towards establishing standards for model evaluation and better reliability of model performance, a collaborative study between two independent teams was created. The Leuven team created one of the most promising femur failure load-predictive model [18], and the Lyon team was in charge of its replication. The objective of this study was to evaluate the reproducibility and replicability of the modeling workflow presented by Sas *et al.* [18], and to study the global sensitivity of this original model.

## 2. Materials and Methods

### *2.1. Datasets*

Two datasets were used. The reproducibility of the workflow was tested using the data presented by the Leuven group [28]. In short, the Leuven dataset [28] consisted of the CT scans (voxel size: 0.4x0.4x0.2 $mm^3$) of 8 femurs with surgical defects (Figure 1) cut at 25 cm (from the top of the head), experimentally tested until failure in a setup mimicking single leg stance. A full description of the experimental setup has been published [29]. A second dataset (the Lyon dataset) was used to test replicability. This dataset consisted of the CT scans (voxel size: 0.7x0.7x1.2 $mm^3$) of 16 intact femurs and 6 femurs (3 pairs) with surgical defects (Figure 1) cut at 15 cm (from the lesser trochanter), together with the description of the experimental setup and experimental failure load used.

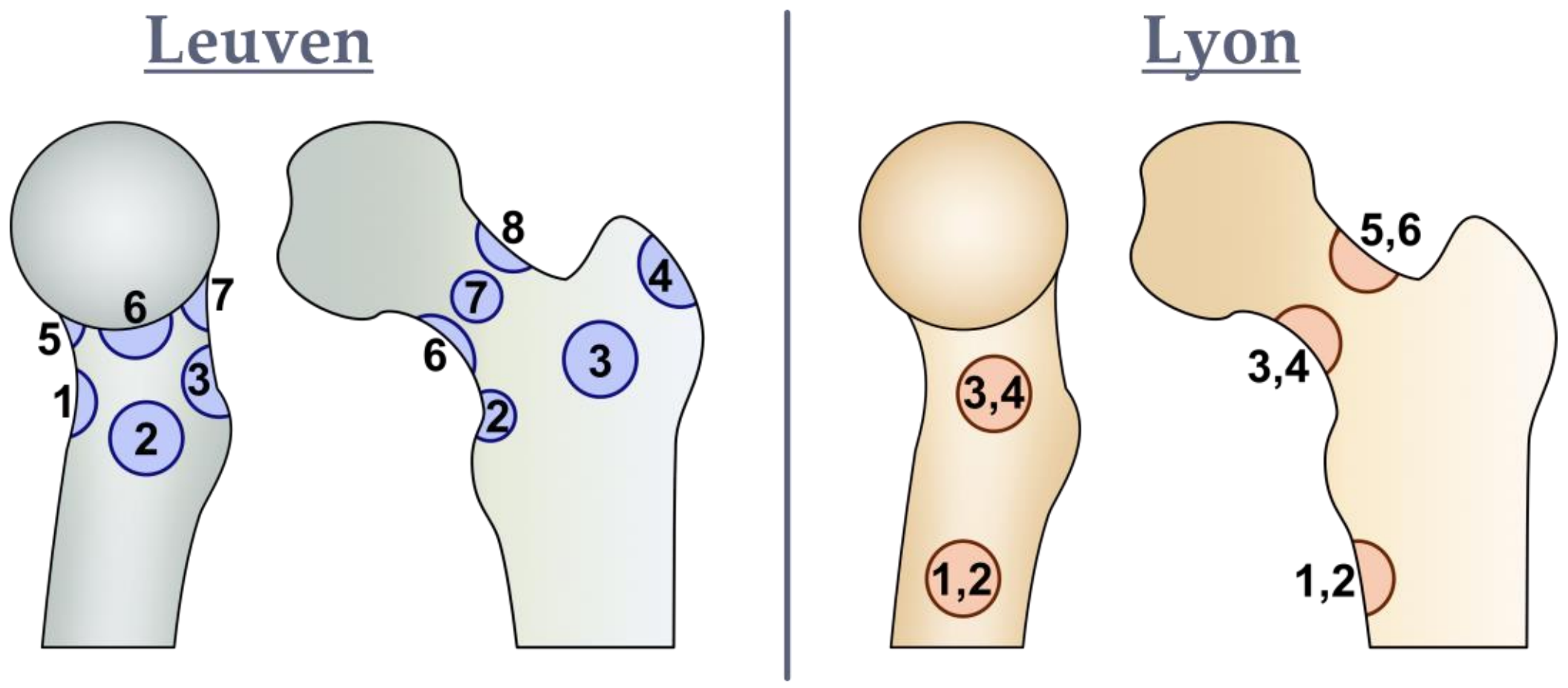

**Figure 1.** Location of defects for both Leuven's and Lyon's samples. Each femur had one defect, the numbers on each defect correspond to the femur it was applied to.

### *2.2. Evaluated workflow*

The evaluated workflow (the Leuven model) is a nonlinear finite element model described and validated in Sas *et al.* 2020 [18]. The model is originally described as follows (*cf.* Figure 2). The femur is segmented from the CT scan (voxel size: 0.94x0.94x3 $mm^3$) (no guideline is given on the segmentation technic). Then, the CT volume is resliced in order to get 2 mm isotropic voxels. A first-order hexahedral (brick, voxel-based) mesh is created based on the segmentation, and mechanical properties are assigned from the Hounsfield unit of each voxel to the associated element. The Hounsfield units are first converted to calcium hydroxyapatite (CaHA) density using a phantom of CaHA scanned along with the femur. Any negative density voxel was set to 0 $g.cm^{-3}$. Surface elements with a CaHA density smaller than 0.25 $g/cm^3$ were removed to account for partial volume effects of the CT scan. Then, the densities were smoothed with a gaussian (SD = 1.27 with limited kernel of 3x3 elements). The CaHA density was then converted into ash density with the relations provided by Keyak *et al.* 2005 [10] (1).

$$\rho_{ash} = 0.887 \times \rho_{CaHA} + 0.0633 \qquad (1)$$

The ash density was then converted into mechanical properties [10].

**Table 1.** Mechanical properties in function of ash density (from Keyak et al. 2005 [10]).

**TABLE 1. Mechanical Property Relationship**

| Relationship | Type of Bone |
|---|---|
| $E$ (MPa) = $14900\rho_{ash}^{1.86}$ | Trabecular and cortical |
| $S$ (MPa) = $102\rho_{ash}^{1.80}$ | Trabecular and cortical |
| $\varepsilon_{AB}$ (mm/mm) = $0.00189 + 0.0241\rho_{ash}$* | Trabecular |
| $\varepsilon_{AB}$ (mm/mm) = $0.0184 - 0.0100\rho_{ash}$* | Cortical |
| $E_p$ (MPa) = $-2080\rho_{ash}^{1.45}$* | Trabecular |
| $E_p$ (MPa) = $-1000$* | Cortical |
| $\sigma_{min}$ (MPa) = $43.1\rho_{ash}^{1.81}$ | Trabecular and cortical |

*To account for the difference between the size of the finite elements and the size of the specimens in which $\varepsilon_{AB}$ and $E_p$ were originally measured, the values of $\varepsilon_{AB}$ and $E_p$ given above must be modified using the following equations before inclusion in the FE models[22]: trabecular bone, $\varepsilon_{AB}' = \varepsilon_{AB} \times (15/3)$ and $E_p' = (3 E E_p)/[15E - (15 - 3)E_p]$; cortical bone, $\varepsilon_{AB}' = \varepsilon_{AB} \times (5/3)$ and $E_p' = (3 E E_p)/[5E - (5 - 3)E_p]$.

In order to accurately reproduce the experimental setup, boundary conditions were carefully applied to the femur model in this study. Specifically, the load was applied on the top surface of a load cup that was placed on the femoral head, while a null displacement was applied to a support embedding the distal diaphysis. These boundary conditions were chosen to closely mimic the experimental conditions, ensuring that the results obtained from the model were as accurate as possible.

To solve the model, the authors used ParOsol, which is an implicit finite element analysis software. This software can simulate linear models. Additionally, a custom modification was applied to the software, as explained in a previously cited paper, to enable the simulation of nonlinear models [26]. Elemental results were computed using gaussian quadrature, which is a reduced integration method commonly used in finite element analysis.

The model performance was evaluated by comparing the predicted failure load with the experimental *ex vivo* failure load. In this study, for easier comparison, we will use the coefficient of determination (**r²**), the average error between predicted and experimental failure loads (**μ**), the standard deviation of errors between predicted and experimental failure loads (**SD**), and the number of samples (**N**) to evaluate the model performance. The model performance in the original paper was ($N_{orig._intact}$ = 10, $r_{orig._intact}^2$ = 0.90, $\mu_{orig._intact}$ = -8%, $SD_{orig._intact}$ = 13%) for intact femurs and ($N_{orig._lesion}$ = 10, $r_{orig._lesion}^2$ = 0.93, $\mu_{orig._lesion}$ = 21%, $SD_{orig._lesion}$ = 22%) for femurs with surgical lesions [18].

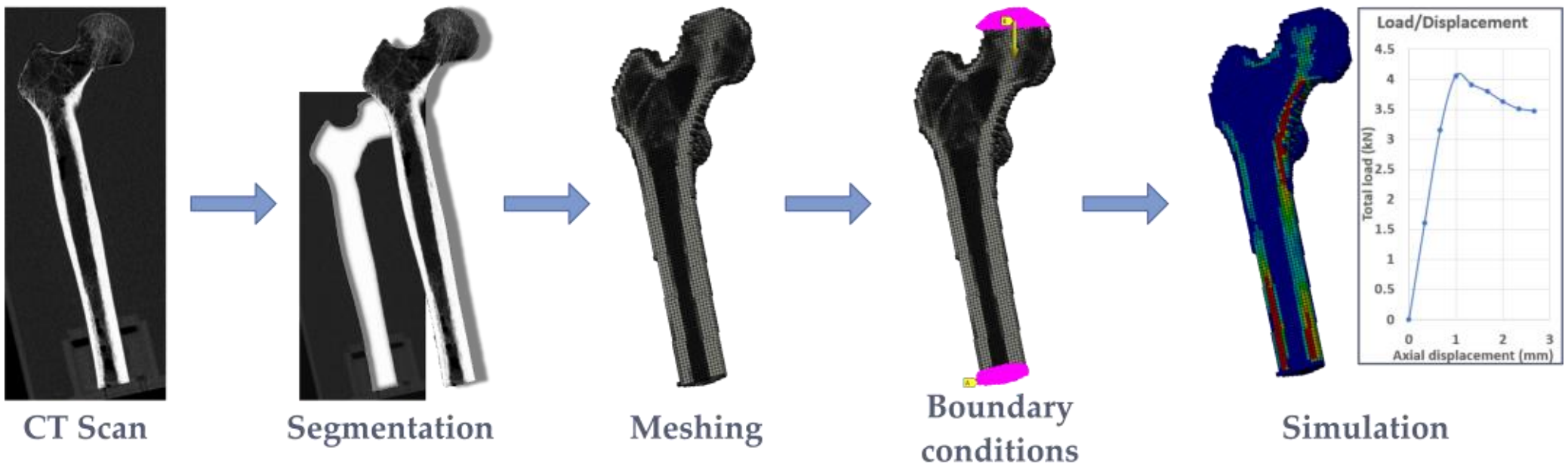

**Figure 2.** Global process to obtain the failure load of a femur based on its CT scan.

Lyon's team first reproduced the model based only on the information provided in the original paper [18]. Leuven's team then reviewed the reproduced model to find any major issue. Some discrepancies between the original and reproduced models are to be noted:

Discrepancies in boundary conditions: the supports were not modelled, the load was directly applied to the top of the femoral head, and the nodes on the surface of the distal diaphysis were applied as a null displacement. While this may have introduced some variability in the results, it has been previously verified that these boundary conditions did not induce any variation on the failure load by modelling the two boundary conditions on three femurs and observing the same failure load for both boundary conditions.

Discrepancies in segmentation technics (no guideline given in the original paper)

Discrepancies in the orientation method of the femurs (Differences between Leuven's team and Lyon's team's orientation method led to less than 2° of difference in the femur orientation)

Discrepancies in the software used. Whereas the original model used Mimics (Materialise NV) and ParOSol [26], the reproduced model used 3DSlicer (slicer.org) for the segmentation, QCTMA (pypi.org/project/qctma) for the mechanical properties attribution, and Ansys Mechanical (ANSYS, Inc.) for the simulation.

These discrepancies were not corrected, as the segmentation and orientation processes were not fully disclosed in the original paper. One major discrepancy was corrected: the original paper's dataset used a CaHA phantom, while the Leuven and Lyon datasets used a $K_2HPO_4$ phantom. Lyon's team used a conversion law to directly get Ash density from $K_2HPO_4$ density [10], while Leuven's team used a conversion law to get CaHA density from $K_2HPO_4$ density [27] and then continued with the originally described process. The reproduced model was then corrected upon choosing the conversion method from Leuven's team.

Each sample of Leuven's and Lyon's datasets was mechanically tested in a compression test (*cf.* **Error! Reference source not found.**Figure 3) using an Instron 3360 (Norwood, MA, USA) for the Leuven's dataset, and an Instron 8802 (High Wycombe, England) for Lyon's dataset. The sample was oriented in a single-leg stance configuration. The sample was preloaded, then small load cycles were applied to remove the nonlinear effects of the initial loading and finally loaded in a quasi-static movement until failure (*cf.* Appendix A for further details). The distal diaphysis was embedded in polymethylmethacrylate (PMMA) and fixed to the testing machine. A PMMA support was molded to each femoral head and placed between the loading plate and the head to distribute the applied load. For the Lyon's dataset, a laser scan was used on each sample to get the exact orientation during the experimental tests and minimize the uncertainty of orientation (the surface mesh of the experimental setup was matched to the surface mesh of the segmentation by a point-set registration algorithm: the coherent point drift using Bayes' inference [30], *cf.* Appendix B) for further details).

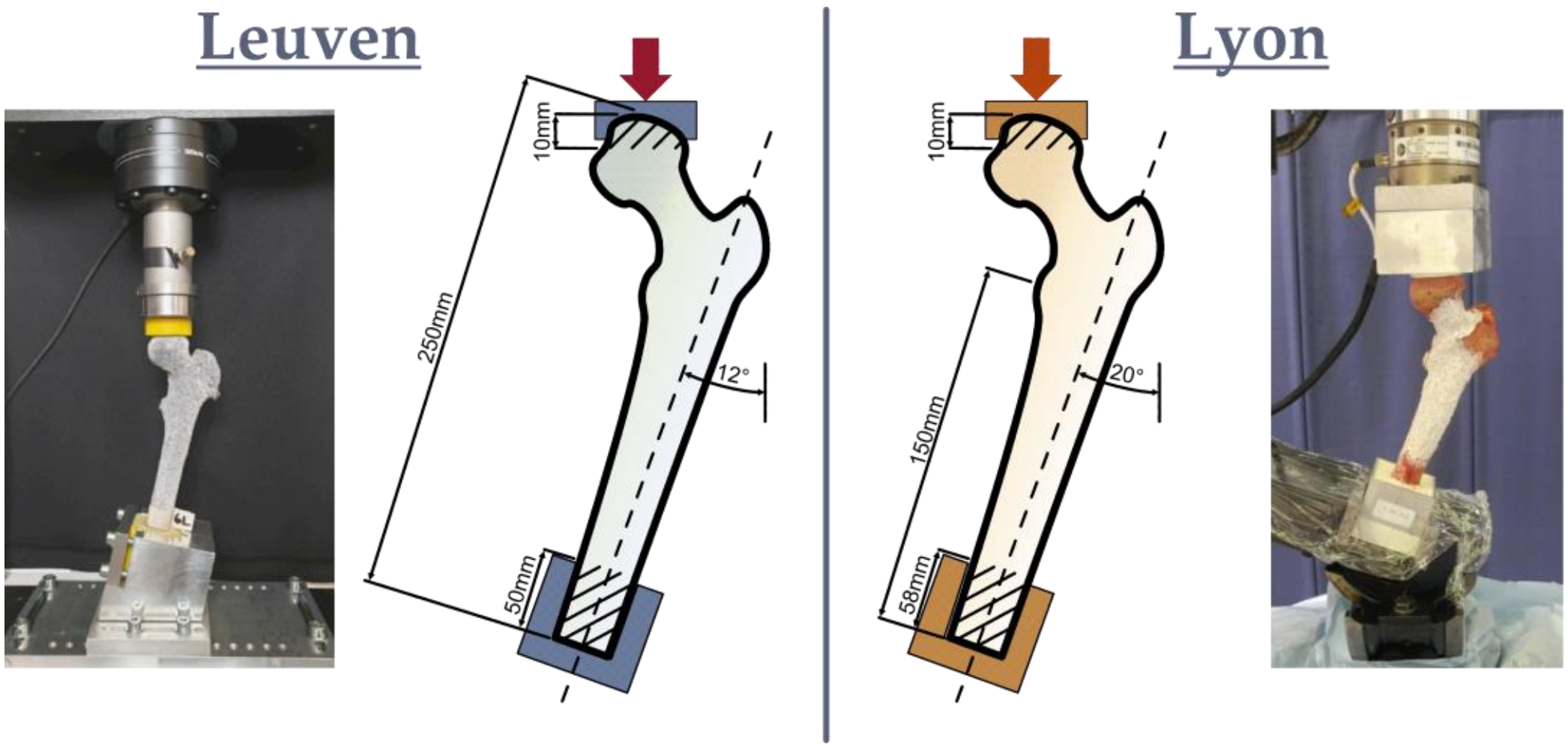

**Figure 3.** Experimental setup description for both Leuven's and Lyon's samples.

Lyon's team applied the reproduced model to each sample of the Leuven dataset, and predictions of the original and reproduced models were compared. To evaluate the replicability of the model, the Lyon's dataset was used. The performance of the reproduced model applied to the Leuven's dataset was then compared to its performance when applied to the Lyon dataset.

To evaluate the performance of models, the mean and standard deviation of the differences between predicted and experimental failure load were used, as were linear regression parameters between predicted and experimental failure loads. The similarity between the original and reproduced model was assessed by the linear regression parameters between the predictions of both models. The value of the differences will be given as a percentage of the relative difference compared to the experimental values.

Two global sensitivity analysis were done to evaluate the reproduced Leuven model. The objective of the first sensitivity analysis was to evaluate the sensitivity of the predicted failure load with regards to the direct input parameters of the model. The mechanical parameters used to define the non-linear constitutional law are obtained with functions dependent of the mineral density of the bone. In that regard, the objective of the second sensitivity analysis was to evaluate the sensitivity of the predicted failure load with regards to the uncertainties of the conversion laws used to convert the mineral density to the mechanical parameters.

The first global sensitivity of the reproduced model was evaluated using the Morris method [35] with 10 trajectories and 4 levels for each of the following parameters: femur length, orientation, segmentation, calibration coefficient and diaphysis rotation (Table 2). The femur used for the global sensitivity analysis was chosen among the Leuven dataset. Ranges of variation for each parameter were chosen by repeatedly applying the nominal values and then computing the standard variation of each measured parameters. The range corresponds to the nominal value +/- two standard deviations. Thus, the results of the Morris method account mostly for the uncertainty due to the operator. The elementary contribution of each parameter on the variation of the failure load will be presented as a percentage of the failure load obtained with the nominal values.

**Table 2.** Parameters studied with the Morris method.

| Parameters | Nominal value | Min | Max | Parameter description |
|---|---|---|---|---|
| **Femur length (mm)** | 250 | 230 | 250 | Working length of the femur |
| **Orientation (°)** | 12 | 6 | 18 | Orientation of the diaphysis axis in the frontal plane |
| **Segmentation** | 0 | -1 | +2 | 0 : reference segmentation<br>-1 : erosion with 1 voxel kernel from the reference<br>+1 an +2 : dilation with respectively 1 and 2 voxels kernels from the reference |
| **Calibration coefficient ($g.cm^{-3}$)** | $8.57 \times 10^{-4}$ | $8.47 \times 10^{-4}$ | $8.67 \times 10^{-4}$ | Slope of the Hounsfield units to $K_2HPO_4$ conversion equation |
| **Diaphysis rotation (°)** | 0 | -8 | 8 | Rotation around the diaphyseal axis |

The second global sensitivity of the reproduced model was evaluated using the Morris method [35] with 5 trajectories and 4 levels for each parameter described in Table 3, totaling 70 simulations. The femur used for the global sensitivity analysis was the same as the one used for the first sensitivity analysis. Ranges of variation for each parameter were chosen by computing the uncertainty of the coefficient of regressions with respect to the method described in Taylor 1997 (Chapter 8) [39]. The uncertainty corresponds to one standard deviation. The range for each parameter used in the sensitivity analysis corresponds to the nominal value +/- the uncertainty.

The sensitivity analyses were conducted by using the SALib Python package (https://github.com/salib/salib).

**Table 3.** Conversion laws of mechanical parameters. Their nominal value and uncertainty were computed from raw data extracted from the article where they were taken from.

| Conversion law | Equation | Parameter | Nominal value | Incertainty | Reference |
|---|---|---|---|---|---|
| $\varrho_{Ash}$ -> E | E = **c_Ee** * $\varrho_{Ash}$ ^ **d_Ee** | c_Ee | 14900 | 320 | [36], [37] |
| | | d_Ee | 1.86 | 0.135 | [36], [37] |
| $\varrho_{Ash}$ -> S | S = **c_S** * $\varrho_{Ash}$ ^ **d_S** | c_S | 102 | 4.4 | [36], [37] |
| | | d_S | 1.8 | 0.201 | [36], [37] |
| $\varrho_{Ash}$ -> $\varepsilon_{AB}$ | $\varepsilon_{AB}$ = **c_epsAB_cort** * $\varrho_{Ash}$ ^ **d_epsAB_cort** | c_epsAB_cort | -0.01 | 0.02998 | [38] |
| | | d_epsAB_cort | 0.0184 | 0.003544 | [38] |
| | $\varepsilon_{AB}$ = **c_epsAB_trab** * $\varrho_{Ash}$ ^ **d_epsAB_trab** | c_epsAB_trab | 0.0241 | 0.02998 | [38] |
| | | d_epsAB_trab | 0.00189 | 0.003544 | [38] |
| $\varrho_{Ash}$ -> $E_p$ | $E_p$ = **Epp_cort** | Epp_cort | -1000 | 1304 | [38] |
| | $E_p$ = **c_Epp_trab** * $\varrho_{Ash}$ ^ **d_Epp_trab** | c_Epp_trab | -2080 | 1304 | [38] |
| | | d_Epp_trab | 1.45 | 0.34 | [38] |
| $\varrho_{Ash}$ -> $\sigma_{min}$ | $\sigma_{min}$ = **c_sigma_min** * $\varrho_{Ash}$ ^ **d_sigma_min** | c_sigma_min | 43.1 | 41.43 | [38] |
| | | d_sigma_min | 1.81 | 0.54 | [38] |

**$\varrho_{K2HPO4}$ : K2HPO4 density**
**$\varrho_{CHA}$ : CHA density**
**$\varrho_{Ash}$ : Ash density**
**E : Young's modulus**
**S : Ultimate stress**
**$\varepsilon_{AB}$ : Max plastic elongation**
**$E_p$ : Softening module**
**$\sigma_{min}$ : Min stress**

## 3. Results

### *3.1. Datasets description*

The ratio of female to male is almost inversed between both datasets with femur with lesions (5:3 for Leuven vs 2:4 for Lyon (lesions)) (*cf.* Figure 4). The Lyon's (intact) dataset has a 7:9 female-to-male ratio. The Lyon (lesions) dataset presents an older population than the other ones, and the femurs show surgical defects at almost half the volume of those of the Leuven's dataset.

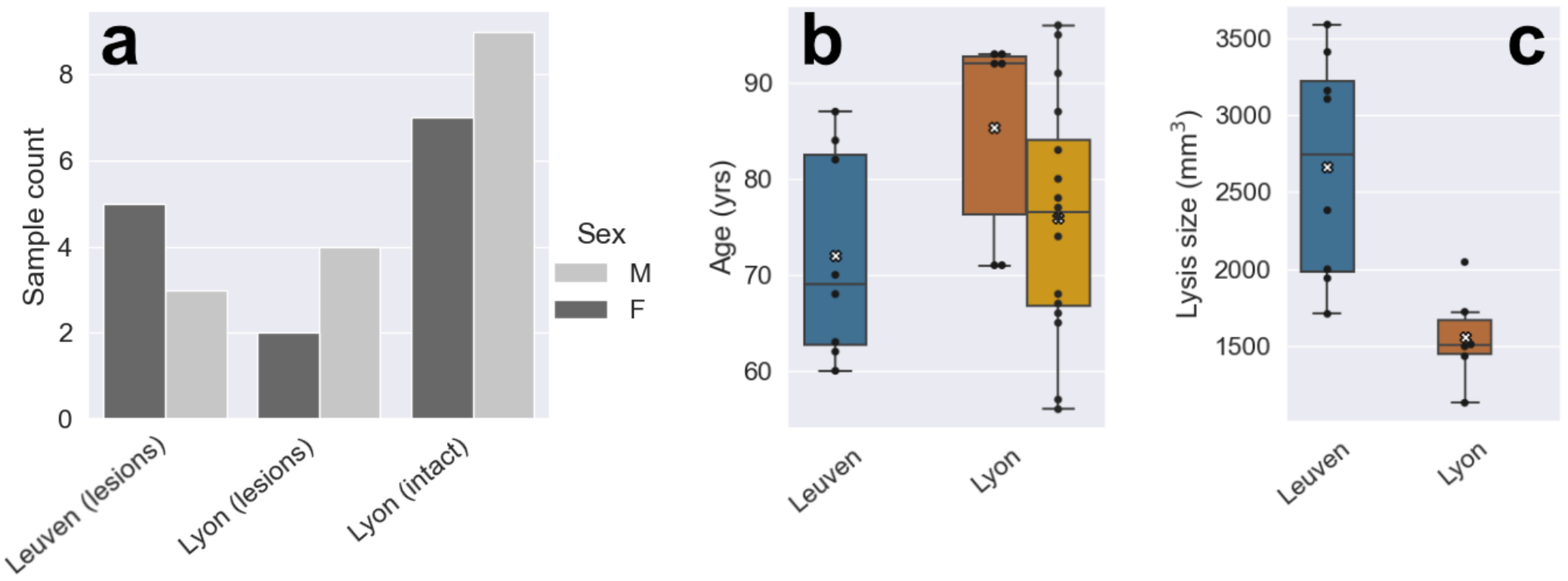

**Figure 4. a**: Sample count differentiated by donor's sex for each dataset, dark grey: female, light grey: male; **b**: Box plot of the donors' age for each dataset; **c**: Box plot of the defect size for datasets with femur with lesions., blue: Leuven's dataset, orange: Lyon's dataset (femurs with lesions), yellow: Lyon's dataset (intact femurs).

### *3.2. Reproducibility*

There is a high correlation between the original and reproduced models ($r_{reproduction}^2$ = 0.95, Figure 5). Dispersions of errors are similar ($SD_{original}$ = 15%, $SD_{reproduced}$ = 16%); however, the systematic bias is doubled for the reproduced model ($\mu_{original}$ = 11%, $\mu_{reproduced}$ = 25%). Correlations with experimental data are similar ($r_{original}^2$ = 0.96, $r_{reproduced}^2$ = 0.98) (Figure 6). Therefore, the model is reproducible, despite having a 14% additional systematic error.

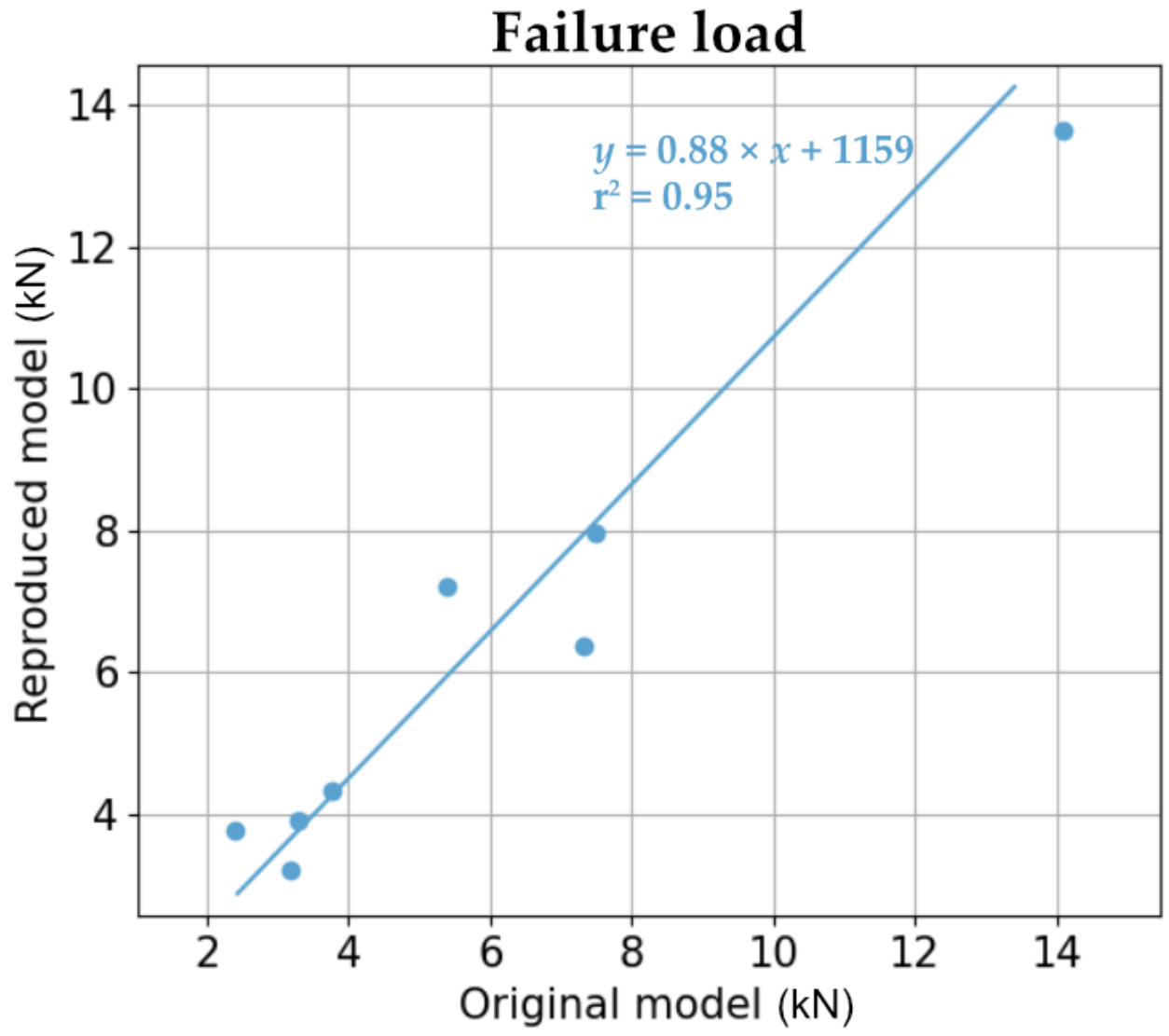

**Figure 5.** Scatter plot and regression of the reproduced vs original model predictions on the Leuven's dataset.

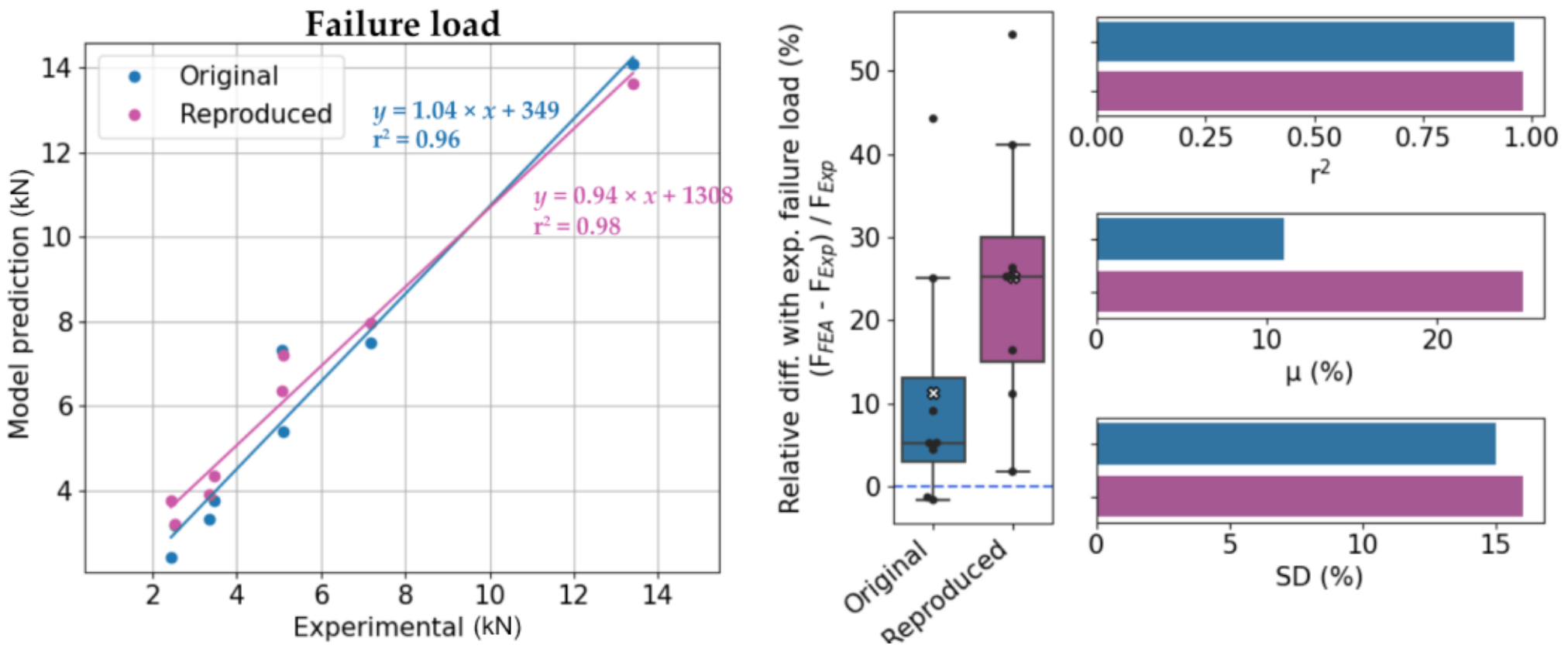

**Figure 6.** From left to right: scatter plot and regression, box plot of the differences, determination coefficient, mean difference and standard deviation of the differences between each model and experimental failure load on the Leuven's dataset. In blue: original model, in purple: reproduced model.

*3.3. Replicability*

The performance of the reproduced model shown previously using the Leuven's dataset tends to degrade when using the Lyon's dataset ($r_{Lyon}^2$ = 0.60, $\mu_{Lyon}$ = -39%, $SD_{Lyon}$ = 28%, Figure 7). These results contain the performance for intact femurs only ($r_{intact}^2$ = 0.66, $\mu_{intact}$ = -50%, $SD_{intact}$ = 17%) and femurs with surgical lesions ($r_{lesion}^2$ = 0.03, $\mu_{lesion}$ = -11%, $SD_{lesion}$ = 32%). There is little to no correlation at all with the experiment for this group. The systematic bias is smaller for femurs with lesions and larger for intact femurs. The dispersion of errors, however, is better for intact femurs and worse for femurs with lesions.

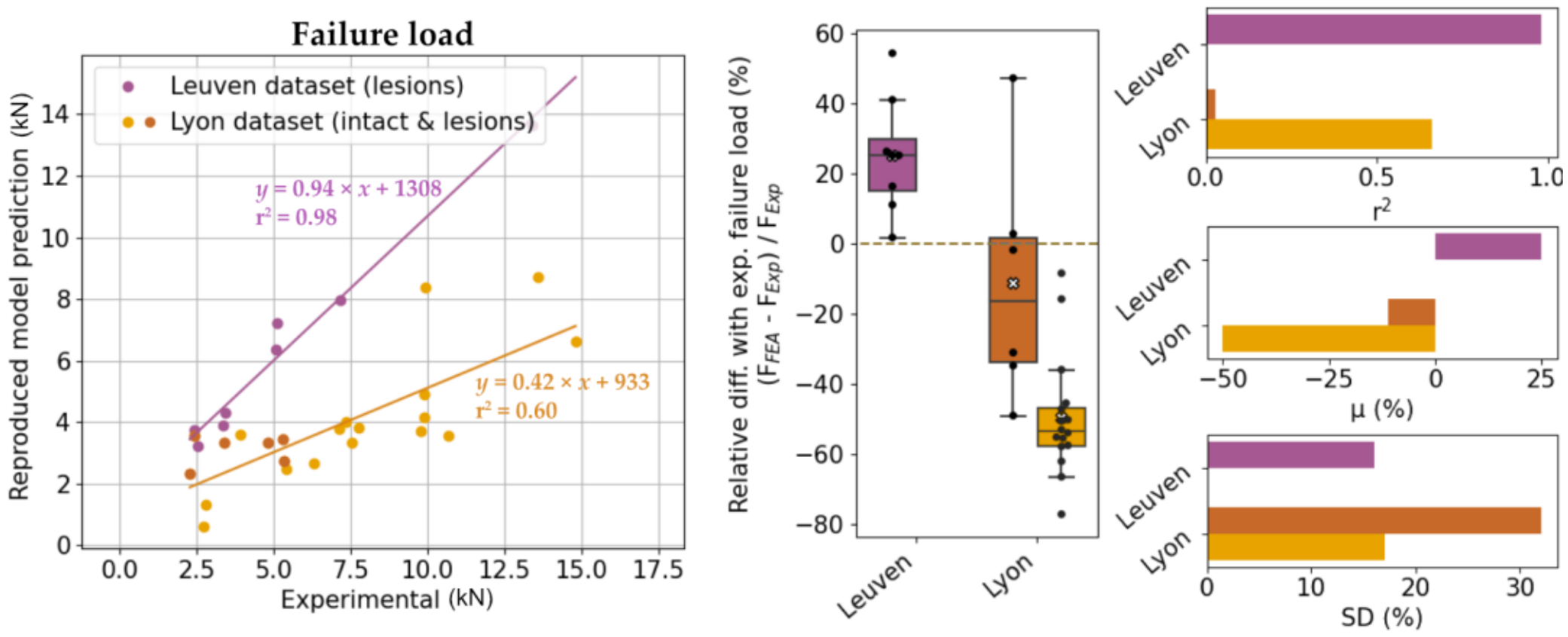

**Figure 7.** From left to right: scatter plot and regression, box plot of the differences, determination coefficient, mean difference and standard deviation of the differences between reproduced model and experimental failure load on each dataset. In purple: Leuven's dataset; in orange: Lyon's dataset (femurs with created lesions); in yellow: Lyon's dataset (intact femurs).

*3.3. Global sensitivity analysis of the input parameters*

Five simulations (spread on 2 trajectories) out of 60 crashed and could not give any result. The 3 defective trajectories were removed from the calculation, leading to the final Morris analysis having 8 trajectories relying on 48 simulations. All the crashed simulations were simulations with the maximum values for the calibration coefficient and the segmentation.

The results are in percentage of the failure load obtained with the nominal values (F = 3766 N). The Morris method shows that the calibration coefficient is the most influent parameter (mean influence of 12.5 % on the failure load) (Figure 8). The orientation has a mean influence of 10.1 %, the segmentation has a mean influence of 8.6 %, the length has a mean influence of 6.6 % and the diaphysis rotation has a mean influence of 4.4 %.

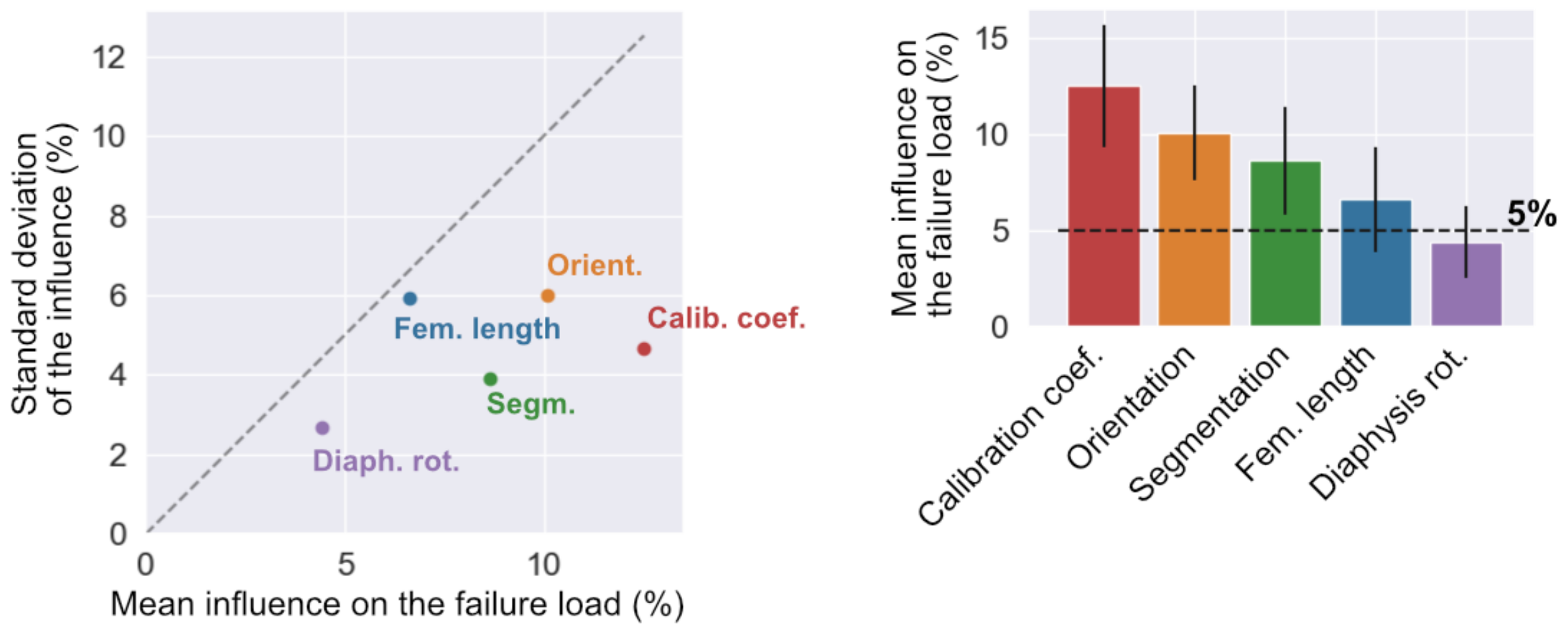

**Figure 8.** Left: Standard representation of the Morris method results (standard deviation of the influence vs mean influence on the failure load, the black dotted line is the y = x curve). The more a parameter is on the right, the more it is influent on the failure load. The more a parameter is high, the more it as a nonlinear influence on the failure load. Right: Ranking of the parameters in function of their mean influence on the failure load. The black vertical line represents the confidence interval of the mean influence computed by bootstrap.

*3.4. Global sensitivity analysis of the conversion laws*

The Moriss analysis of the parameters of the conversion laws shows influences superior to 5% on the failure load of 6 parameters out of 13 (Figure 9). The parameters of the ultimate stress (d_S and c_S) have a mean influence of respectively 44 % and 14 %. The parameters of the minimum stress (d_sigma_min and c_sigma_min) have a mean influence of respectively 26 % and 7 %. The parameter of the softening modulus of the cortical (Epp_cort) has a mean influence of 12 % and finally, one parameter of the plastic max elongation (c_epsAB_cort) has a mean influence of 10 %.

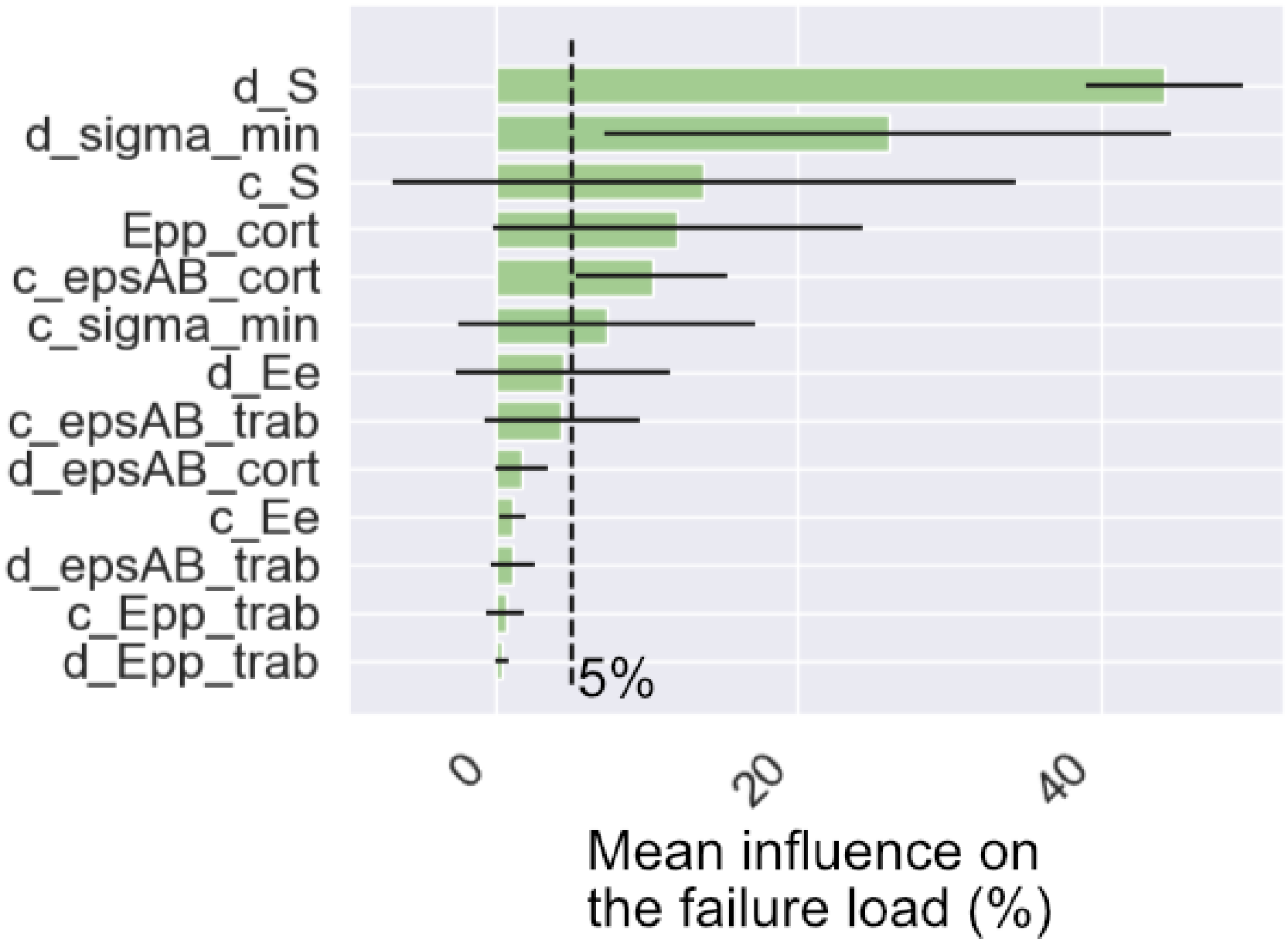

**Figure 9.** Bar plot of the mean influence of each parameter of the conversion laws. The black horizontal lines represent the standard deviation of the mean influence.

## 4. Discussion

### *4.1. Reproducibility*

This study mainly aimed to evaluate the reproducibility and replicability of a femur model predictive of failure load from the literature. As described in Baker's article [22], the chances to perfectly reproduce experiments are low. Despite not giving the exact same predictions, the original model from Leuven's team and the reproduced model from Lyon's team gave very similar results on Leuven's dataset, except that the reproduced model gave systematically higher failure loads (by 14%). The origin of these discrepancies lies in the unknown influence that some parameters might have on the model. The segmentation and orientation methods differed between the original and reproduced models. As the results of the reproducibility analysis suggest, the major difference is the systematic bias, hinting of a systematic variation between both models (*e.g.,* a systematic "over segmentation" or a systematic angle between both models for each sample). Taddei and colleagues [23] showed that, to some extent, there were uncertainties in the predictions of their femur model induced by the geometry of the mesh and the material properties. Despite the plurality of segmentation technics [31], material properties [32] (*i.e.,* density calibration and mechanical properties relations), and femur models in the literature [8 - 18] [31] [33], no proper global sensitivity analysis is required when publishing such models, leading to inevitable reproducibility issues if all parameters are not fully constrained by exhaustive and complete explanation of every step of the model creation.

### *4.2. Replicability*

The replicability study showed a huge difference in performance depending on the dataset. The original model on the Leuven's dataset showed similar performance as with the intact data of the original paper (both performance results were obtained by the original author). However, the results of the reproduced model obtained on the Lyon's dataset gave worse precision and correlations for both intact femurs and femurs with lesions. This could question the quality of the reproduction of the model. As the reproducibility results showed, there is a high correlation between the reproduced and original models on the same dataset (Leuven's), hinting that the performance decrease does not come from the reproduction quality, but from the dataset itself or from the experimental settings. Although similar, there are factors changing between the experimental setup of both datasets which potentially affect model precision: femurs are not the same effective length (the length from the top of the femoral head to the distal diaphysis above the embedding is not the same), and the orientation changes by 8° between both setups. The donors' characteristics and their femur properties also changed slightly (*cf.* Figure 4). The surgical defect size could influence the performance of the model, although the performance was good for both intact and femurs with lesions in the original paper. Thus, the difference of performance might be better explained by other parameters which were not measured and not reflected by the spatial distribution of the mineral density, such as mineral quality or micro porosity.

*4.3. Global sensitivity analysis*

The calibration coefficient is the most influent parameter (12.5 % in average on the failure load). This may be explained by the systematic use of the density to compute every mechanical parameter (*e.g.,* Young's modulus, ultimate stress, etc.) of the nonlinear constitutive law.

Each law defining the conversion between density and mechanical parameters were obtained by regression of experimental data. These regressions intrinsically have uncertainties due to the inability of the density to explain 100% of the mechanical parameter, and the small sample size (N = 10) used for the laws of the maximum plastic elongation, the softening modulus and the minimum stress.

By computing the uncertainty based on the raw data used by Keyak *et al.* 2005 to create the nonlinear law used by the Leuven model, it appears the uncertainty of the law gets high (Figure 10).

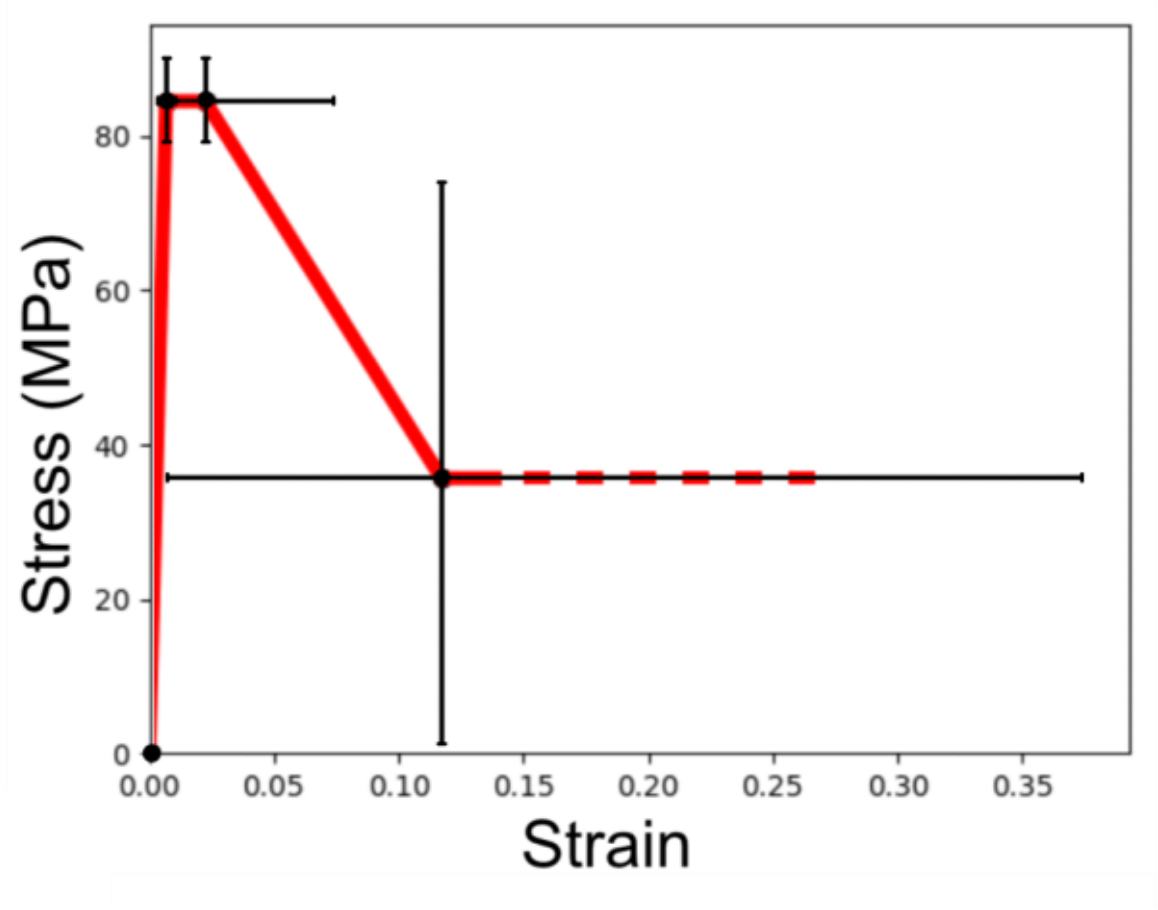

**Figure 10.** Nonlinear constitutive law used in the Leuven model (red) for an ash density of 1.0 g/cm$^3$ (cortical bone). The black lines represent the uncertainty of the position for each point (*cf*. Table 3).

The standard deviations of the mean influence are close to, sometimes higher than, the mean influence. This could be interpreted by a high nonlinearity and interaction between parameters, however the low number of trajectories used for the Morris analysis could explain the high standard deviations. Finally, the laws giving the ultimate stress and the minimal stress use the coefficients among the 6 most influent towards the failure load. In priority, these laws should be better approximated in order to be more confident in the predicted failure load of the model.

The uncertainty of the mechanical parameters highly influences the predicted failure load. The cause of the high influence of the calibration coefficient might come from the uncertainties of the definition the mechanical parameters as the uncertainty of the calibration coefficient is propagated by the laws of conversion of the mechanical parameters. To illustrate this, the Figure 11 shows a sample of the laws defining the Young's modulus found in the literature. The definition of the most basic parameter used for a mechanical model is far from being unanimous and characterizes the high uncertainties of the input parameters used in the mechanical models of bones.

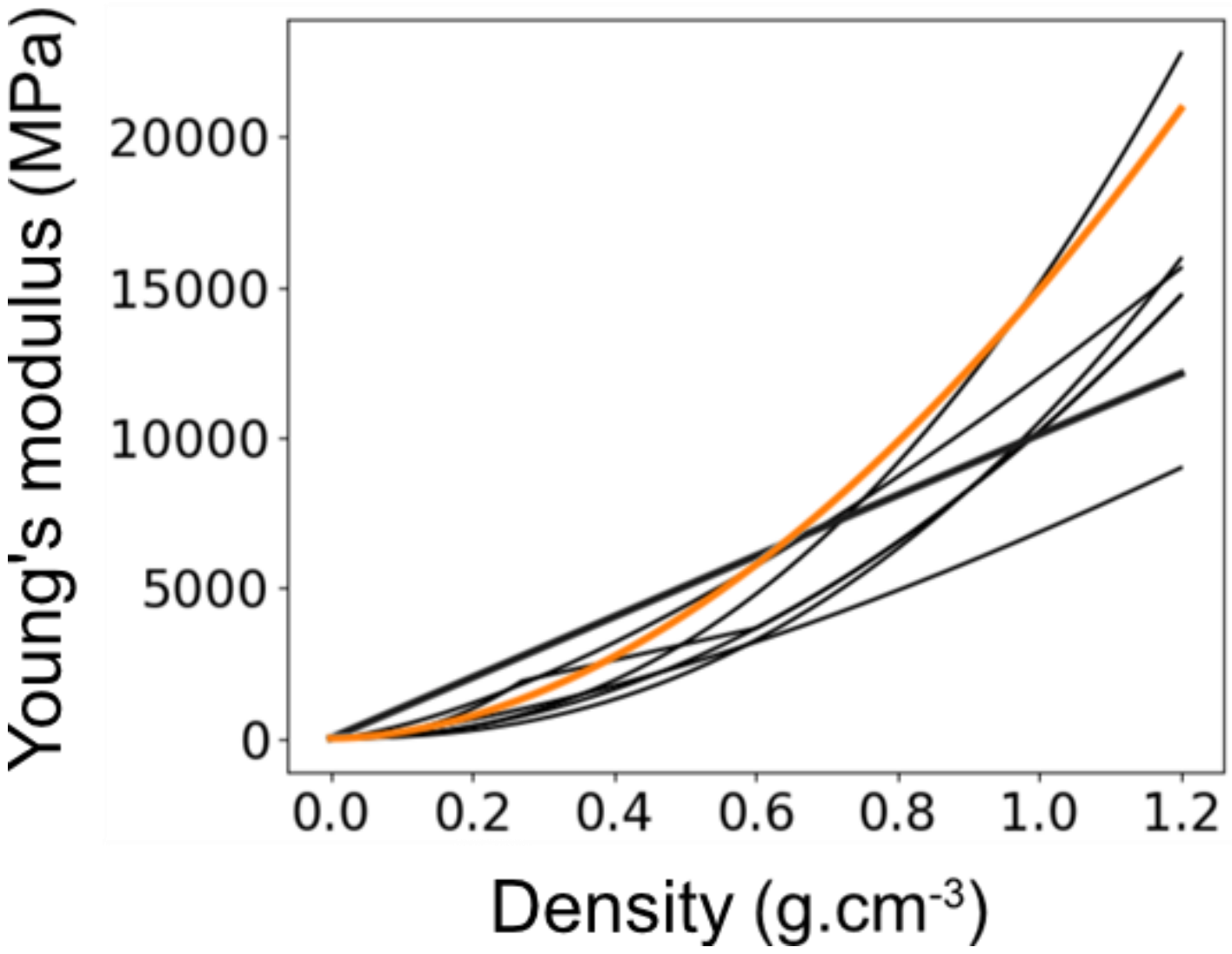

**Figure 11.** Laws of conversion from density to Young's modulus extracted from Helgason *et al.* 2016 [40]. The Orange curve is the law used by the Leuven model.

The limitations of this study are notable and should be considered when interpreting the results. One of the main limitations is the small sample size used in some groups of the analysis (femurs with lesions), which may have affected the reproducibility of the model. A larger sample size would have allowed for more statistically significant results, providing greater confidence in the findings. Additionally, this study did not have the opportunity to evaluate the Lyon samples using the original model, which could have shed light on the robustness of the model and its applicability to other datasets.

Another shortcoming of this study is the use of different FE solvers. Whereas the original model used ParOsol, the reproduced model used Ansys. While unlikely, the use of different software packages could have had an impact on the variability of the failure loads.

It is important to note that the limitations of this study do not invalidate its findings, but rather provide valuable information for future research. By addressing these limitations, future studies could provide more robust and reliable results, ultimately advancing our understanding of the biomechanical processes involved in this field. Finally, the difference in performance depending on the dataset is commonly the cause of overfitting when creating the model. However, the dataset used in the original paper [18] and the Leuven's dataset gave similar performance, which indicates a lesser probability for the overfitting cause.

## 5. Conclusions

This study showed through reproducibility and replicability analysis that, although CT-scan-based finite element models can predict fairly correctly femur failure load, there are still some unknown influences from parameters such as length or orientation, from methods such as segmentation, or from characteristics of the femurs that may alter the model performance. The knowledge of these influences is currently not required to validate a model, thus very few if any of the already existing models study these influences through a global sensitivity analysis coupled with uncertainty quantification. Sensitivity and uncertainty studies are necessary if these models are to be used in clinical applications, as femurs characteristics and imaging modalities may differ. The credibility assessment of such models is a topic of growing concern for the medical and scientific community using them [34]. In order to better understand the behavior of this model, an uncertainty quantification should be done to correctly evaluate the confidence of the predictions. Finally, this study showed that automatic process of the measurement may be needed for the parameters of segmentation, orientation and femur length in order to minimize the influence of the operator on the predictions.

**Author Contributions:** Conceptualization, H. Follet, M. Gardegaront, D. Mitton; methodology, M. Gardegaront, A. Sas; software, M. Gardegaront, A. Sas; validation, F. Bermond, D. Brizard, H. Follet, M. Gardegaront, D. Mitton, A. Sas, H. Van Lenthe; formal analysis, M. Gardegaront; investigation, M. Gardegaront; resources, F. Bermond, D. Brizard, H. Follet, D. Mitton, H. Van Lenthe; data curation, M. Gardegaront; writing—original draft preparation, M. Gardegaront; writing—review and editing, F. Bermond, D. Brizard, C. Confavreux, H. Follet, M. Gardegaront, A. Levillain, D. Mitton, J. B. Pialat, A. Sas, H. van Lenthe; visualization, D. Brizard, M. Gardegaront; supervision, F. Bermond, D. Brizard, C. Confavreux, H. Follet, A. Levillain, D. Mitton, J. B. Pialat, H. van Lenthe; funding acquisition, C. Confavreux, H. Follet, D. Mitton, J.B. Pialat.

All authors have read and agreed to the published version of the manuscript.

**Funding:** This work was partly funded by LabEx Primes (ANR-11-LABX-0063) and MSD Avenir Research Grant.

**Institutional Review Board Statement:** For Leuven femurs, ethical approval was granted by the Ethics Committee of the University Hospitals Leuven (reference number NH019 2018-09-02). The Lyon femurs were provided by the Departement Universitaire d’Anatomie Rockefeller (Lyon, France) through the French program on voluntary corpse donation to science (DC-2015-2357).

**Conflicts of Interest:** The authors declare no conflict of interest.

**Data Availability Statement:** Raw data for the Lyon dataset may be acquired at: https://data.univ-gustave-eiffel.fr/dataverse/ex-vivi-experiment-on-bones

## Appendices

**Appendix A**. Details of CT-scan and mechanical testing for each dataset.

| | **Leuven** | **Lyon** |
|---|---|---|
| *CT-Scan* | • Scanned in water basin<br>• Scanned with $K_2HPO_4$ phantom (Mindway Model 3)<br>• Scanned with distal support<br>• Siemens Somatom Force, Siemens AG, Germany<br>• 120 ref kV, 250 ref mAs<br>• Slice thickness 0.4 mm, slice increment 0.2 mm, pitch 0.85, in-plane resolution 0.4 mm and bone kernel | • Scanned without water basin<br>• Scanned with $K_2HPO_4$ phantom (Mindway Model 3)<br>• Scanned before embedding<br>• 120 kV, 200 mAs<br>• Slice thickness 1.25mm, pitch: 1, in-plane resolution 0.74 |
| *Mechanical testing* | • Preload of 50 N<br>• 20 sinusoidal preconditioning cycles (50-500 N, 1Hz)<br>• Load rate of 10 N/s | • Preload of 20 N<br>• 4 preconditioning cycles of 0.4 mm loading plate displacement<br>• Displacement rate of 0.21 mm/s |

**Appendix B**. Surface matching process.

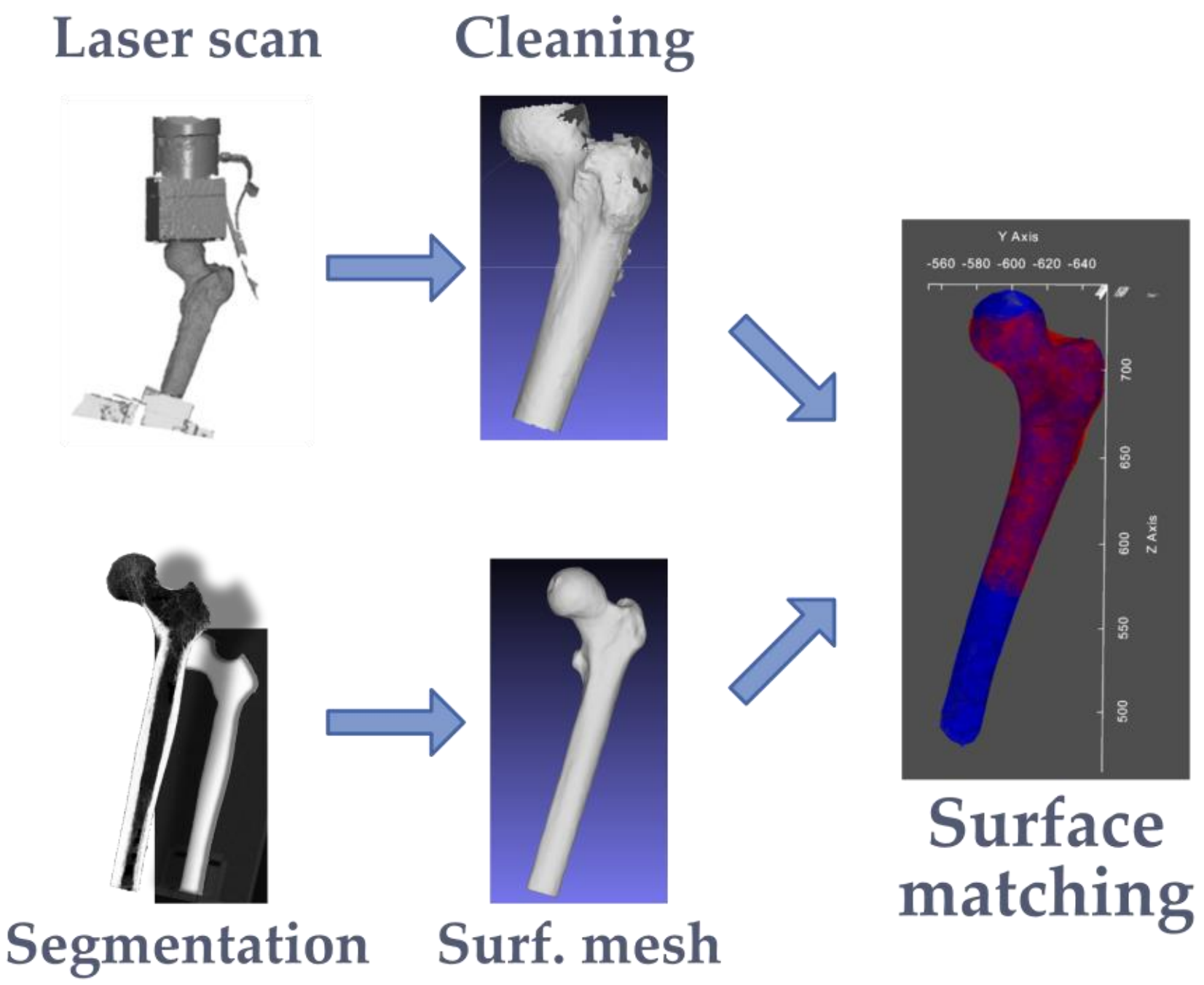

**Declaration of interests**

☐The authors declare that they have no known competing financial interests or personal relationships that could have appeared to influence the work reported in this paper.

☒The authors declare the following financial interests/personal relationships which may be considered as potential competing interests:

| |
|---|
| Marc GARDEGARONT reports financial support was provided by LabEx PRIMES. Marc GARDEGARONT reports financial support was provided by MSD Avenir. If there are other authors, they declare that they have no known competing financial interests or personal relationships that could have appeared to influence the work reported in this paper. |